\documentclass[12pt,twoside]{article}
\usepackage{amsmath,amsfonts,amssymb,theorem}
 \date{9 Novembre  1998}
\title{ Un th\'eor\`eme de Nakai-Moishezon pour certaines
classes de type (1,1) }
\author{Philippe Eyssidieux}

  \newtheorem{theo}{Th\'eor\`eme}
 \newtheorem{prop}{Proposition}
 \newtheorem{coro}[prop]{Corollaire}
 \newtheorem{lem}[prop]{Lemme}

 \theorembodyfont{\rmfamily}
 \newtheorem{defi}[prop]{D\'efinition}

\newenvironment{prv}{\paragraph{Preuve}}{\par\medskip}

 \newcommand{\Si}{\Sigma}

\begin{document}
 \maketitle
\begin{abstract} Let $X$ be a smooth compact projective variety over $\mathbb C$.

Let $H^2(\pi_1(X),\mathbb R)^{1,1}$ be the intersection of
$H^{1,1}(X,{\mathbb R})$ with the image of the map $H^2(\pi_1(X),{\mathbb R})\to H^2(X)$
induced by the classifying map $X\to B\pi_1(X)$.
Let $NS(X)$ be the N\'eron-Severi group of $X$.

Let $[\omega]\in H^2(\pi_1(X),\mathbb R)^{1,1}+ NS(X)\otimes {\mathbb R}$.
In this note, we prove that  $[\omega]$ is the cohomology class of a  K\"ahler
metric  if and only if
for every $d$-dimensional reduced closed algebraic subvariety  $Z\subset X$,
$[\omega]^d.Z>0$.
\end{abstract}

Un r\'esultat classique de  Harvey et Lawson  \cite{HL} affirme que,
si $(M,\omega_M)$ est une vari\'et\'e K\"ahlerienne connexe
avec $\dim(M)=d$,
une classe
$\alpha\in H^{1,1}(M)$ est une classe de K\"ahler ssi
$\alpha.\omega_M^{d-1}>0$ et pour tout $(d-1,d-1)$-courant positif $dd^c$-ferm\'e
$T$ $\alpha.T>0$.

Les cycles utilis\'es pour tester l'amplitude sont juste $dd^c$-ferm\'es, ce qui contraste avec
le classique th\'eor\`eme de Nakai-Moishezon et son extension aux diviseurs r\'eels \cite{CP}.
Le r\'esultat de \cite{CP} est le suivant:
Soit $X$ une  variet\'e  compacte projective sur $\mathbb C$.
Soit $NS(X)$ le groupe de N\'eron-Severi de $X$.
 toute classe $[\omega] \in NS(X)\otimes {\mathbb R}$
v\'erifiant $[\omega]^d.Z>0$ pour tout sous espace analytique
 $Z$  r\'eduit de dimension $d$ de $X$
est une classe de K\"ahler.

Le contraste entre ces deux r\'esultats est frappant et laisse ouverte la possibilit\'e
de th\'eor\`emes de type Nakai-Moishezon pour des classes de cohomologie r\'eelles arbitraires.
La question qui motive cette recherche est la suivante.
Est-il vrai que toute classe $[\omega]$ de type $(1,1)$ sur  une vari\'et\'e projective alg\'ebrique $X$
 v\'erifiant
 $[\omega]^d.Z>0$
 pour tout sous espace analytique $Z$  de dimension $d$ de $X$
est une classe de K\"ahler?

Un  r\'esultat r\'ecent de A. Lamari -et ind\'ependamment N. Buchdahl-, annonc\'e dans \cite{L}
et bas\'e sur une exploitation astucieuse de \cite{HL}, l'affirme quand $X$ est une surface.

Soit $X$ une vari\'et\'e projective alg\'ebrique complexe. On note
 $H^2(\pi_1(X),\mathbb R)^{1,1}$ l'intersection de $H^{1,1}(X,\mathbb R)$
avec l'image de l'application $H^2(\pi_1(X),\mathbb R) \to H^2(X,\mathbb R)$
induite par l'application classifiante $X\to B\pi_1(X)$.

D\'efinissons
$\overline{ H^2(\pi_1(X),\mathbb R)^{1,1}+  {NS}(X)\otimes {\mathbb
R} }$ comme le $\mathbb R$-vectoriel form\'e des classes de type
$(1,1)$ $\omega$  sur  $X$ telles qu'il existe une famille
projective lisse sur un polydisque  $\Xi\to \Delta$ avec $\Xi_0=X$,
 une famille continue $(\omega_t)_{t\in \Delta}$ de classes de
type $(1,1)$ avec  $\omega_0=\omega$ et  une suite $t_n$ avec $t_n
\to 0$ telles que
  $\omega_{t_n} \in NS(\Xi_{t_n})\otimes {\mathbb R}+H^2(\pi_1(\Xi_{t_n}),\mathbb R)^{1,1}$.

\begin{theo}\label{theo1}

Soit $X$ une vari\'et\'e complexe projective lisse compacte.

Soit $[\omega]\in \overline{ H^2(\pi_1(X),\mathbb R)^{1,1}+  {NS}(X)\otimes {\mathbb R} }$.
 $[\omega]$ est la classe de cohomologie d'une m\'etrique de K\"ahler
ssi pour chaque sous espace alg\'ebrique  $d$-dimensionnel r\'eduit  $Z\subset X$,
$[\omega]^d.Z>0$.
\end{theo}

Les techniques standard de recollement et r\'egularisation de courants
rappel\'ees dans la section 1.3 r\'eduisent, par r\'ecurrence
sur la dimension de $X$,  la  question \`a montrer que la classe
$[\omega]$ est big et v\'erifie le lemme de Kodaira,  c'est \`a dire
peut \^etre repr\'esent\'ee par
 un courant strictement positif ferm\'e \`a singularit\'es logarithmiques.

 L'hypoth\`ese $[\omega]\in H^2(\pi_1(X),\mathbb R)^{1,1}+ NS(X)\otimes {\mathbb R}$
signifie que la classe $[\omega]$ est limite de la premi\`ere classe de Chern   d'une suite de fibr\'es
lin\'eaires holomorphes
d\'efinis sur le rev\^etement universel de $X$ et sur lesquels agit une extension idoine de
$\pi_1(X)$ par $S^1$.
Que $[\omega]$ est big s'obtient  de fa\c con tr\`es classique en utilisant
une adaptation  d\'evelopp\'ee dans \cite{E2}
du th\'eor\`eme d'indice $L_2$ d'Atiyah \cite{Ati} qui donne des techniques
 cohomologiques d'\'etude des groupes de sections $L_2$ de tels fibr\'es.

Un exemple d\^u \`a R. Livne \cite{Liv}, montre qu'il existe une
surface projective $S$ uniformis\'ee par l'espace hyperbolique complexe (i.e.:
$c^2_1(S)=3c_2(S)$) avec $\overline{NS(S)\otimes \mathbb R}=NS(S)\otimes \mathbb R \not = H^{1,1}(S)$ (la premi\`ere
\'egalit\'e vient du fait que $S$ est rigide.). Pourtant notre
caract\'erisation du c\^one K\"ahlerien s'applique \`a $S$
puisque c'est un $K(\pi,1)$, ce qui implique que $H^{1,1}(S,{\mathbb R}) =
H^2(\pi_1(S),\mathbb R)^{1,1}$.

Je remercie M. Paun pour d'utiles discussions
concernant les techniques de r\'egularisation et recollement
de courants.

\section{Courants quasipositifs et m\'etriques singuli\`eres pour les espace complexes r\'eduits}

\subsection{Courants sur un espace complexe r\'eduit}

\paragraph{}
Soit $(S,O_S)$ un espace complexe r\'eduit.  Narasimhan d\'efinit le faisceau
$\mathcal{PSH} \cap C^0$ des fonctions plurisousharmoniques continues
de $(S,O_S)$, comme suit: une fonction continue est plurisousharmonique s'il existe
un recouvrement de $S$ par des ouverts de cartes $(U_i)_i$  munis de plongements
$U_i\subset \mathbb C^{N_i}$ tels que $\phi_{U_i}$ est la restriction \`a $U_i$
d'un fonction plurisousharmonique continue sur un voisinage de $U_i$ dans $\mathbb C^{N_i}$. Que cette d\'efinition est
consistante est d\'emontr\'e dans \cite{Nara}.

De la m\^eme fa\c con, on d\'efinit le faisceau  $C^{\infty}_{\mathbb R}$ des fonctions lisses r\'eelles,
le faisceau $\mathcal H _{\mathbb R}$ des fonctions pluriharmoniques r\'eelles, le faisceau $\mathcal{PSH}$
des fonctions  plurisousharmoniques (on demande que $\phi$ soit semi continue sup\'erieurement et que
pour toute composante irr\'eductible locale $T$ $\phi_{T}\not\equiv -\infty$), le faisceau des fonctions psh lisses $\mathcal{PSH} \cap C^{\infty}$ ,
le faisceau $\mathcal{QPSH}=\mathcal{PSH}+C^{\infty}_{\mathbb R}$  des fonctions quasi psh .

\cite{Nara} d\'efinit aussi le faisceau $\mathcal{SPSH}$ des fonctions strictement plurisousharmoniques
comme \'etant le faisceau des fonctions $\phi$ telles pour toute fonction continue $f$
il existe $\epsilon >0$ tel que $\phi+\epsilon f$ soit psh. Les fonctions strictement psh lisses
sont localement les restrictions par un plongement de fonctions lisses dont la forme de Levi est
 sup\'erieure \`a une forme de K\"ahler.

\paragraph{}
Soit $F$ un faisceau en $\mathbb R$-vectoriels  et $G$ un sous faisceau en ${\mathbb R}$-vectoriels
Soit $K\subset F$ un sous faisceau d'ensembles, $G$-invariant. On note $K/G$  l'image faisceautique
de $K$  dans
$F/G$.

Une {\em m\'etrique K\"ahl\'erienne} sur $(S,O_S)$  est une section globale de  $\mathcal{SPSH}\cap C^{\infty}_{\mathbb R}/{\mathcal H}_{\mathbb R}$.
Un {\em courant lisse} est une section globale de $C^{\infty}_{\mathbb R}/ {\mathcal H}_{\mathbb R}$.
Un {\em courant big} est une section globale de $\mathcal{SPSH}/\mathcal H _{\mathbb R}$.
Un {\em courant psef} est une section globale de $\mathcal{PSH}/ \mathcal H _{\mathbb R}$.
Un {\em courant quasi-positif} est une section globale de $\mathcal{QPSH}/ \mathcal H _{\mathbb R}$.

 Soient $T$ un courant quasi-positif  et $\gamma$ un courant lisse. Par d\'efinition,
$T\ge \gamma$ ssi $T-\gamma$ est psef.

Par analogie avec les notations habituelles, on notera par $dd^c$ l'application naturelle
$C^{\infty} _{\mathbb R} \to C^{\infty} _{\mathbb R} / \mathcal H _{\mathbb R} $.

\subsection{M\'etriques singuli\`eres}

\paragraph{}
On a des applications  $H^0(S,\Phi /{\mathcal H}_{\mathbb R})\to
 H^1(S,\mathcal H_{\mathbb R})$ surjectives quand $\Phi= C^{\infty},\mathcal{QPSH}$,
d'images le {\em c\^one K\"ahler} quand $\Phi= \mathcal{SPSH}\cap C^{\infty}$,
le {\em c\^one big }quand $\Phi=\mathcal{SPSH}$.

Soit $L$ un fibr\'e holomorphe en droites. $O_S(L)$ est un faisceau
inversible. Une {\em m\'etrique hermitienne singuli\`ere} sur $L$
est une application $h$ de l'espace total $L$ dans $[0,\infty[$ qui
sur chaque fibre est une m\'etrique hermitienne possiblement nulle
s'\'ecrivant localement $h= e^{-\phi}$ avec $\phi$ quasi psh. Soit
$h$ une m\'etrique singuli\`ere sur $L$. La collection des
potentiels locaux $\phi$ de $h$ d\'efinit un courant quasi positif
qu'on est en droit d'appeler $C_1(L,h)$. La classe de $C_1(L,h)$
dans $H^1(S,\mathcal H_{\mathbb R})$ ne d\'epend que de $L$.

Si $S$ est compact $H^1(S,\mathcal H_{\mathbb R})$ est un $\mathbb R$-vectoriel de dimension finie.

\subsection{Lemme de recollement}

\paragraph{}
Soit $(S,O_S)$ un espace complexe r\'eduit compact et $\omega_S$ une forme
 de K\"ahler sur $S$ \footnote{On peut formuler une variante sans supposer que  $S$ est un espace
K\"ahl\'erien.}.

Le lemme de r\'egularisation de Richberg \cite{R} et le lemme  de recollement de Paun \cite{P}
sont tous les deux valides pour les espaces complexes r\'eduits, par la  preuve originelle:

\begin{lem} Soit $T$ un courant quasipositif sur $S$, section globale de $\mathcal{QPSH}\cap C^0/{\mathcal H}_{\mathbb R}$
et $\gamma$ un courant lisse avec $T\ge \gamma$.

Pour tout $\epsilon>0$,  il existe un repr\'esentant lisse $T_{\epsilon}$ de $[T]$ avec $T_{\epsilon }\ge \gamma -\epsilon \omega_S$.
\end{lem}

\begin{lem}\label{glue}
Soit $Z\subset S$ un sous espace complexe compact. Soient $\alpha$, $\gamma$ deux courants lisses
sur $S$.

Soit $\phi_1$ une fonction quasi psh sur $S$ lisse hors de $Z$
et $\phi_2$ une fonction lisse dans un ouvert $U$ contenant $Z$.
On suppose $\alpha +dd^c\phi_1\ge \gamma$ sur $S$
et $\alpha|_U +dd^c\phi_2 \ge \gamma|_U$.

Alors pour tout $\epsilon>0$, il existe $\phi_{\epsilon}\in C^{\infty }(S)$  tel que
$\alpha+dd^c\phi_{\epsilon} \ge \gamma -\epsilon \omega_S$.

\end{lem}

Le  lemme suivant est une cons\'equence directe de la preuve de \cite{D1}, Theorem 4, p. 285
(voir aussi \cite{P}, lemme 1).

\begin{lem}
\label{prolposi}
Soit $Z\subset S$ un sous espace complexe compact. Soient $\alpha$, $\gamma$ deux courants lisses
sur $S$. On suppose que $[\alpha|_Z]\in H^1(Z,{\mathcal H}_{\mathbb R})$ poss\'ede un repr\'esentant lisse
$\alpha_Z$ tel que $\alpha_Z\ge \gamma|_Z$. Alors, pour tout
$\epsilon >0$,  il existe un voisinage ouvert $U$ de $Z$
dans $S$ et $\phi_U\in C^{\infty}(U)$ tel que $\alpha|_U+dd^c\phi_{U,\epsilon}
\ge \gamma|_U -\epsilon \omega_S|_U$.
\end{lem}

\begin{prv} Puisque \cite{D2} et \cite{P} ne formulent pas les choses de cette fa\c con,  recopions l'argument.

On peut supposer que $\alpha|_Z=\alpha_Z$.

Consid\'erons une chaine de sous espaces analytiques compacts  de $S$
$\emptyset =Z_{-1} \subset Z_0 \subset Z_1\subset \ldots\subset Z_N=Z$
avec $ Sing(Z_i) \subset Z_{i-1}$.

Supposons, par r\'ecurrence, construits $\phi_i$ une fonction lisse sur $S$
 et $V_i$ voisinage de $Z_i$ tels que:
$$(\alpha +dd^c\phi _i )|_{V_i}
\ge(\gamma -\frac{\epsilon}{ 2^{k-i}} \omega_S)|_{V_i}$$
$$(\alpha +dd^c\phi _i )|_{Z}
\ge(\gamma -\frac{\epsilon}{ 2^{k-i}} \omega_S)|_{Z}$$

Soient $U$ un ouvert de $S$ ne recontrant pas $Z_{i+1}$ et
 $(V_{\lambda},U_{\lambda})$ une famille finie de cardinal $\Lambda$ d'ouverts de $S$ telle que
$V_{\lambda} \Subset U_{\lambda}$, $ \mathfrak V=(V_{\lambda})_{\lambda}\cup V_i\cup U$ est un recouvrement ouvert de $S$ et il existe des plongements $i_{\lambda}:U_{\lambda} \to \mathbb  C^{N_{\lambda}} $
 avec
$$Z_{i+1}\cap U_{\lambda}=i_{\lambda}^{-1} ( \{z_{1}=0, \ldots, z_{N_{\lambda}-d}=0\})$$
Soit $f_{\lambda}$ une fonction lisse sur $S$ \`a valeurs dans $C^{N_{\lambda}-d}$
coincidant avec
$ (i_{\lambda}^*z_k)_{k=1}^{k=N_{\lambda}-d}$
sur $V_{\lambda}$. Soit $((\theta_{\lambda})_{\lambda}, \theta_i,\theta'_i)$ une partition de l'unit\'e lisse
subordonn\'ee \`a $\mathfrak V$.

Pour $1\gg \eta >0$, d\'etermin\'e durant la preuve, on pose:

 $$\phi_{i+1}= \phi_i + \sum_{\lambda} \theta_{\lambda}   \eta^3 \log(1+\eta^{-4}|f_{\lambda}|^2)  $$

$dd^c( \theta_{\lambda}\log(1+\eta^{-4}|f_{\lambda}|^2) ) $
est somme de quatre termes.
Le premier terme $dd^c( \theta_{\lambda})\log(1+\eta^{-4}|f_{\lambda}|^2)$
est born\'e par $M.\omega_S$, le second
$d\theta_{\lambda}\wedge 2(d^c f_{\lambda},f_{\lambda})/(1+\eta^{-4}
|f_{\lambda}|^2)$
est born\'e par $M\eta^{-2}\omega_S$ ($M$ est une constante positive assez grande.). Le troisi\`eme terme
est similaire au second.

Le quatri\`eme terme
$\theta_{\lambda}dd^c\log(1+\eta^{-4}|f_{\lambda}|^2)$
est positif.

Quitte \`a diminuer $\eta$,
 sur  $W_i=\{\theta_i>1/3\}$,  on a bien $\alpha+dd^c\phi_{i+1} \ge \gamma -\epsilon/ 2^{k-1-i}\omega_S$
 et $(\alpha+dd^c\phi_{i+1})_Z\ge (\gamma -\epsilon/ 2^{k-1-i}\omega_S)_Z$.

Par ailleurs si $s\in Z_{i+1}\cap
\{ \theta_i\le 1/2 \}$ on peut choisir $\lambda$ tel que $\theta _{\lambda}>1/3\Lambda$
pr\`es de $s$. Par hypoth\`ese on a, pr\`es de $s$,  pour $C$ une constante positive:

$$\alpha+ dd^c\phi_i -\gamma \ge-\frac{\epsilon}{2^{k-i}}\omega_S -C( i df_{\lambda}\wedge d\bar f_{\lambda})$$

Quitte \`a diminuer encore $\eta$ on peut supposer que, dans l'ouvert $V_s$
contenant $s$ :
$$ \frac {1}{3\Lambda} \eta^{3}dd^c \log(1+\eta^{-4}|f_{\lambda}|^2)
= \eta^{-1} \frac{ ( i df_{\lambda}\wedge d\bar f_{\lambda})}
 {3\Lambda (1+\eta^{-4} |f_{\lambda}|^2)^2}
 \ge
C ( i df_{\lambda}\wedge d\bar f_{\lambda})$$

On peut choisir $s_1,\ldots, s_M$ tel que $(V_{s_i})_i$ est un recouvrement ouvert
de $Z_{i+1}\cap \{\theta_i \le  1/2\}$  et par suite
$\eta$ tel que sur $V_{i+1}= \{\theta_i >1/3\} \cup V_{s_1} \cup\ldots\cup V_{s_M}$:
$$
(\alpha +dd^c\phi_{i+1})_{V_{i+1}} \ge( \gamma-\frac{\epsilon}{2^{k-1-i}}\omega_S )|_{V_{i+1}}
$$
Ceci conclut la preuve du lemme \ref{prolposi}.
\end{prv}

 Le lemme suivant est d'une utilit\'e \'evidente dans les questions de type
Nakai-Moishezon:

\begin{lem}\label{recoll} Soit $(S,O_S)$  un espace complexe compact K\"ahlerien.
Soit $[\omega]\in H^1(S,{\mathcal H}_{\mathbb R})$. On suppose que $[\omega]$
est repr\'esent\'ee par un courant big qui est  lisse en dehors d'un ensemble analytique propre
$E$
et que $[\omega|_E]\in H^1(E,{\mathcal H}_{\mathbb R})$ est dans le cone K\"ahler de $E$.
Alors $[\omega]$ est dans le c\^one K\"ahler de $S$.
\end{lem}
\begin{prv}
Cons\'equence imm\'ediate des lemmes \ref{glue} et \ref{prolposi} .
 \end{prv}

Parmi les techniques de r\'egularisation et de recollement de courants quasi positifs d\'evelopp\'ees par Demailly \cite{D2},\cite{D3},\cite{P}, seul le difficile  th\'eor\`eme de r\'egularisation \cite{D2} et ses cons\'equences
ne s'\'etendent pas
imm\'ediatement au cas d'une base singuli\`ere. Il est
n\'eanmoins naturel de conjecturer une version singuli\`ere de
\cite{D2}.

\section{Cohomologie $l_2$ des faisceaux p\'eriodiques}

\subsection{Faisceaux analytiques coh\'erents p\'eriodiques
sur un rev\^etement galoisien infini d'une vari\'et\'e alg\'ebrique}

\paragraph{}
Soit $\Sigma$ la donn\'ee d'un groupe discret  $\Gamma$ op\'erant proprement discontinuement par biholomorphismes sur un espace complexe analytique $(S,O_S)$. $\Si_{red}$
d\'esignera l'espace complexe r\'eduit associ\'e sur lequel $\Gamma$ agit  par biholomorphismes.
Soit $\{1\} \to S^1 \to G\to \Gamma\to \{1 \}$ une extension centrale de $\Gamma$ par $S^1$.

Un {\em faisceau analytique coh\'erent $G$-p\'eriodique} (en abr\'eg\'e $G$-fac) sur le $G$-espace $(S,O_S)$ est un faisceau analytique coh\'erent sur $(S,O_S)$  muni d'un rel\'evement de l'action naturelle de $G$.
Un {\em morphisme de $G$-fac} est un morphisme de faisceaux $O_S$-lin\'eaire commutant \`a $G$.

Soit $\chi$ un caract\`ere continu  de $S^1$ .
Un {\em $G,\chi$-fac} est un $G$-fac $F$ tel que,  pour tout point $s$ de $S$,  l'action de $S^1$ sur l'espace $F_s$ des germes
en $s$ de sections de $F$ soit donn\'ee par la caract\`ere $\chi$. La sous cat\'egorie pleine de
la cat\'egorie des $G$-fac dont les objets sont les $G,\chi$-fac  est une cat\'egorie
ab\'elienne $C_{G,\chi}(\Sigma)$.

\paragraph{}
Soit $X$ une vari\'et\'e k\"ahlerienne et $\pi:\tilde X \to X$ son
rev\^etement universel. Soit $ [\omega]\in H^2(\pi_1(X),\mathbb
R)^{1,1}$ et $\omega $ une forme lisse ferm\'ee de bidegr\'e $1,1$
repr\'esentant cette classe de cohomologie. Il existe alors une
extension centrale $\{ 1 \} \to S^1 \to G \to \pi_1(X)\to \{1\}$ du
groupe fondamental de $X$ par $S^1$, un fibr\'e lin\'eaire
holomorphe $(\tilde L,\bar\partial,h)$ sur $\tilde X$ et un
rel\'evement \`a $(L,\bar\partial,h)$  de l'action de $G$ sur
$\tilde X$ induite par l'action naturelle de $\pi_1(X)$ tels  que
la premi\`ere forme de Chern-Weil de  $(L,\bar\partial,h)$ soit
$\omega$ \cite{G}, \cite{E0} . $S^1$ agit par le caract\`ere $\chi$.   Le
faisceau des sections $\mathcal L = O_{\tilde X} (\tilde L)$ est un
$G,\chi$-fac sur $\tilde X$. Le foncteur $\pi^* . \otimes \mathcal
L^m$ donne une \'equivalence de cat\'egories de la cat\'egorie de
faisceaux analytiques coh\'erents sur $X$ vers $C_{G,\chi^m}
(\tilde X)$ .

\subsection{Foncteurs cohomologiques}

\paragraph{}
On suppose  $S/{\Gamma}$ compact.
Pour tout $p\in [1,\infty]$, il est possible de d\'efinir l'espace $H^0_{p}(\Sigma,\mathcal F)$ des sections $L^p$
d'un objet $\mathcal F$ de $C_{G,\chi}(\Sigma)$, \cite{E1}.

Il est \'egalement possible  de prolonger cette d\'efinition en
 construisant des groupes de
cohomologie $L_p$
$H^q_{p}(\Sigma,{\mathcal F})$, nuls si $q \not \in \{ 0,\ldots, \dim S
\}$.

Ces groupes de cohomologie s'organisent en
un $\delta$-foncteur cohomologique (voir \cite{Har} III.1, p.205), ce qui signifie
qu'\`a toute suite exacte courte de $C_{G,\chi}(\Sigma)$
correspond une suite exacte longue de cohomologie $L_p$.

\subsection{Th\'eorie d'indice $L_2$ d'Atiyah}

$H^q_{p}(\Sigma,{\mathcal F})$ est  un $G$-module topologique, non
n\'ecessairement Hausdorff si $q\not=0$. Quand $p=2$, on peut de
plus estimer la taille de ces groupes de cohomologie. En effet, le
plus grand quotient Hausdorff de $H^q_{2}(\Sigma,{\mathcal F})$, $\bar H^q_{2}(\Sigma,{\mathcal
F})$,
est un objet d'une sous  cat\'egorie additive tr\`es particuli\`ere de la cat\'egorie des
 repr\'esentations hilbertiennes de $G$, sur laquelle existe une fonction
 dimension $\dim_G$,  v\'erifiant les propri\'et\'es
usuelles de la fonction dimension en alg\'ebre lin\'eaire classique au d\'etail
pr\`es qu'elle est \`a valeurs r\'eelles.
\footnote{La d\'ecouverte de cette fonction $\dim_G$ est d\^ue \`a Murray et Von
Neumann.}

Soit $f: \Sigma ' \to \Sigma$ un morphisme propre $\Gamma$-\'equivariant et
 $\mathcal F \in Ob C_{\bar G}(\Sigma ')$. On dispose d'une suitre
 spectrale de Leray-Serre qui donne lieu \`a la relation de d\'evissage:
 $$
 \sum_q (-1)^q \dim_G \bar H^q_2(\Sigma',\mathcal F) =\sum_{p,q}(-1)^{p+q} \dim_G \bar H^p_2(\Sigma,R^pf_*\mathcal F)
 $$

Cette relation donne lieu \`a une proc\'edure de calcul de
l'invariant $\chi_2(\Sigma, \mathcal F)=\sum_{i=0}^{\dim S} \dim_G
\bar H^q_{2}(\Sigma,\mathcal F)$ par r\'ecurrence sur la dimension
$d$
du support de $\mathcal F\in Ob C_{\bar G}(\Sigma)$.

 Un d\'evissage simple ram\'ene au cas
o\`u $\Sigma$ est r\'eduit.  Par
\cite{Hiro} et \cite{Rossi}, il existe une d\'esingularisation
\'equivariante $\mu:\Sigma'\to \Sigma$ telle que
 $\mathcal V= \mu^*{\cal F}/ T$ soit localement libre o\`u $T$ est
  le sous faisceau de torsion maximal de $\mu^*{\cal F}$.

Il suit qu'il existe deux familles finies $((S^{\pm}_0,\mathcal V^{\pm}_0), \ldots, (S_d^{\pm},\mathcal V_d^{\pm}))$
o\`u $S_i^{\pm}$ est une $\Gamma$-vari\'et\'e
propre cocompacte
avec
$\dim S_i^{\pm}=i$ et $\mathcal V^{\pm}_i $ un faisceau $G$-\'equivariant
localement libre sur $S_i^{\pm}$ telles que :

$$ \chi_2 (\Sigma,\mathcal F)=\chi_2(S_d^+,\mathcal V_d^+)+ \sum_{q=0}^{d-1} \chi_2(S_i^{+},\mathcal
V^{+}_i) -\chi_2(S_i^{-},\mathcal
V^{-}_i)
$$

  Dans le cas o\`u le faisceau \'equivariant  $\mathcal V$ est
  localement libre, si $\chi$ est le caract\`ere trivial, le r\'esultat de \cite{Ati}
  permet de calculer $\chi_2(\Sigma,\mathcal V)$. M\^eme si $\chi$ n'est pas le
  caract\`ere trivial,
  un argument alternatif invoquant le  th\'eor\`eme d'indice local de Getzler implique (voir \cite{E0} pp. 194-195):

$$\chi_2(\Sigma', \mathcal V) =\int _{S'/\Gamma} \text{ch}(\mathcal V) \text{Todd}  (T_{S'})$$

Pour tout ceci, voir \cite{E2} (La th\'eorie d'indice $L_2$ d'Atiyah  \cite{Ati} ne
suffit pas. L'adaptation de \cite{Ati} dans  \cite{E1}
est insuffisante puisqu'elle se limite au cas o\`u $S$ est projective  lisse.).

Il est observ\'e dans \cite{E1} que les preuves analytiques des
th\'eor\`emes d'annulation usuels de la g\'eom\'etrie alg\'ebrique
complexe (Kodaira-Akizuki-Nakano, Grauert-Riemenschneider,
Kawamata-Viehweg, ...) fonctionnent dans le cas de la cohomologie
$L_2$.

\section{Amplitude d'un faisceau inversible projectivement p\'eriodique}
\paragraph{}
Soit $L_{G,\chi}(S)$  la sous cat\'egorie pleine de $C_{G,\chi}(\Sigma)$ form\'ee
des objets dont le faisceau analytique coh\'erent
sous jacent est inversible.  On note $L_G (\Sigma)=\cup_{\chi \in X(S^1)}L_{G, \chi}$,
$C_G(\Sigma)=\cup_{\chi \in X(S^1)}L_{G, \chi}$...

\begin{defi}

Soit $(S,O_S)$ un espace complexe r\'eduit. $\mathcal L$ est dit $G$-nef si et seulement si
la classe de cohomologie de ${\mathcal L}$ dans $H^1(S,{\mathcal H}_{\mathbb R})$ est dans l'adh\'erence du c\^one K\"ahler.

En g\'en\'eral $\mathcal L$ est dit $G$-nef sur $\sigma$ ssi $\mathcal L_{red}$
est $G$-nef sur $\Sigma_{red}$ .
\end{defi}

\begin{defi}
\label{defiample} On suppose $\chi\in \{ 0, \pm Id \}$. ${\mathcal L}\in Ob (L_G(\Sigma))$ est dit $G$-ample  si et seulement si
et pour tout
pour tous $\mathcal F \in Ob (C_G(\Si))$,  $\mathcal N \in N^+_G(\Si)$ et $q>0$ il existe $N_{\cal F}\in
\mathbb Z$ tel que si $n\ge N_{\cal F}$ on ait
$H^q_{2}(S, \mathcal F \otimes {\mathcal N} \otimes \mathcal L ^{\otimes n})=0$
\end{defi}

\begin{lem}
${\mathcal L}\in Ob (L_G(\Sigma))$  est $G$-ample ssi $\mathcal L_{red}\in Ob (L_G(\Sigma_{red})$  est $G$-ample.
\end{lem}

\begin{lem} Si $\mathcal L$ est  un faisceau inversible $G$-ample, il existe
 une m\'etrique hermitienne lisse  $G$-invariante
sur $\mathcal L$ dont la forme de courbure est une forme k\"ahlerienne.
\end{lem}

\begin{lem}Si $\mathcal L$ est  $G$-ample, $\mathcal L$ est  $G$-nef.
\end{lem}

\begin{lem}  \label{quasitriv} Soit $\Sigma$ est un rev\^etement galoisien infini d'un espace complexe
projectif $(S,O_S)$. On suppose que $ Ob (L_{G,id})\not =\emptyset$.
Soit  $L$ un fibr\'e lin\'eaire  ample (resp. nef) sur $F$.
Alors,
$\pi^*O_S(L)$ est $G$-ample (resp.$G$-nef).
\end{lem}

\begin{prv} Nous  pouvons supposer, par r\'ecurrence sur $d=\dim_{\mathbb C}
 \Sigma$,
que, pour tout $G,\chi$-fac $\mathcal F$  dont le support est de
dimension $< d$,
il existe $n_{\mathcal F}$ tel que, si  $n\ge n_{\mathcal F}$ et
$q\ge 1$,  $H^q_2(\Sigma, \mathcal F  \otimes \mathcal N
\otimes \pi^*O(nL))=0$ , pour  tout $G,\chi$-faisceau inversible $G$-nef
$\mathcal N$.

Soit $\omega^{GR}_S$  le faisceau canonique de Grauert-Riemenschneider
de $S$.
Soit $\mathcal W$ un $G,\chi$-fac de support $\Sigma$ de rang $r$.  Soit $H$ un
diviseur de Cartier ample sur $S$ et ${\mathcal L}' \in Ob L_{G,\chi}$
tels qu'il existe $ (\omega^{GR}_S (-H))^{\oplus r}\otimes {\mathcal L}' \to
\mathcal W $ un morphisme qui est g\'en\'eriquement un
isomorphisme.  Cette fl\'eche donne lieu \`a deux suites exactes,

$$0\to \mathcal K \to (\omega^{GR}_S (-H))^{\oplus r}\otimes {\mathcal L}' \to \mathcal J\to 0$$
$$0\to \mathcal J \to \mathcal W \to \mathcal C\to 0$$

$\mathcal C$ et $\mathcal K$ ont des supports de dimension $<d$.

On invoque le th\'eor\`eme d'annulation de Grauert-Riemenschneider
pour d\'eduire que, pour $n\ge n_{H} $ et $q\ge 1$,
 $H^q_2(\Sigma, \pi^* \omega^{GR}_S (-H)\otimes\mathcal N \otimes \pi^*O_S(nL) )=0$.

La suite exacte longue de cohomologie et l'hypoth\`ese de
r\'ecurrence pour $\mathcal K$ transf\'erent l'annulation
asymptotique
\`a $\mathcal J$. Le m\^eme argument donne  que, pour $n\ge n_{\mathcal W} $ et $q\ge 1$,
$H^q_2(\Sigma, \mathcal W\otimes \mathcal N \otimes
\pi^*O_S(nL))=0$. Ceci conclut la preuve.

\end{prv}

\section{Preuve du Th\'eor\`eme \ref{theo1} }

\subsection{Cas o\`u $[\omega]\in H^2(\pi_1(X),\mathbb R)^{1,1}+ NS(X)\otimes {\mathbb Q}$. }
\paragraph{}
On suppose que $\pi:\Sigma\to S$ est le rev\^etement universel d'un espace complexe projectif alg\'ebrique
compact $S$ de dimension $d$.

\begin{prop}Soit ${\cal L} \in Ob(L_{G,\chi}(\Si)$.
Soit $H$  une section hyperplane de $S$. On pose $\tilde H= H\times _S\Si$.

Supposons que  $\mathcal L_{H\times _{S}\Si}$
est un objet $G$-ample de $L_{G,\chi}(\tilde H)$.

Pour tous $\mathcal F \in Ob (C_G(\Si))$,  $\mathcal M \in N^+_G(\Si)$ et $q>1$ il existe $N_{\cal F}\in
\mathbb Z$ tel que si $n\ge N_{\cal F}$ on ait
$H^q_{2}(S, \mathcal F \otimes {\mathcal M} \otimes \mathcal L ^{\otimes n})=0$.
\end{prop}

\begin{prv} Fixons $m_0\in \mathbb Z$.
Soit $s$ la section tautologique de $O_S(H)$. On a la suite exacte dans $C_{G,\chi}(\Sigma)$:
$$ 0\to \mathcal F\otimes \pi^*O_S(-H) \buildrel{\pi^*s}\over{\to}\ \mathcal F \to \mathcal F_H \to 0$$
${\mathcal F}_H=i_*\Phi$ o\`u $\Phi$ est un objet de $C_{G}(\tilde H)$.
Plus g\'en\'eralement on a la suite exacte courte obtenue en tensorisant  par
$ \pi^*O_S((m+1)H)\otimes \mathcal N
\otimes \mathcal L ^n
$. Utilisant la suite exacte longue de cohomologie $L_2$ associ\'ee \cite{E2} et
la d\'efinition \ref{defiample} avec $\mathcal N=\mathcal M_H\otimes \pi^*O_H(m H)$,
il existe
$N_{\mathcal F}$ tel que  si $n\ge N(m_0,\mathcal F)$, $m\ge 0$ et $q\ge
2$:
$$H^q_{2} ( \Si,{\mathcal F}\otimes {\mathcal M}\otimes \pi^*O_S(mH)
 \otimes {\mathcal L}^n)\simeq H^q_{2} ( \Si , {\mathcal F}\otimes {\mathcal M}\otimes \pi^*O_S((m+1)H) \otimes {\mathcal L}^n) $$
Faisant tendre $m$ vers l'infini, en utilisant le lemme \ref{quasitriv} , il suit que
 $H_2^q(\Sigma, {\mathcal F}\otimes {\mathcal M}\otimes {\mathcal L}^n)=0$
pour $q\ge 2$ et $n\ge N_{\mathcal F}$.
\end{prv}

Utilisant le th\'eor\`eme de Riemann-Roch pour la cohomologie $L_2$ \cite{E2}, suit le:

\begin{coro} Si, de plus ${\mathcal L}^d.S >0$, pour tout
$\mathcal F$ de $C_{G}(\Sigma)$ qui n'est pas de torsion,
il existe $c>0$ tel que:

$$\dim_G H^0_2(\Sigma,\mathcal F \otimes {\mathcal L}^n )= c n^{d} +O(n^{d -1})$$
\end{coro}

Ici, on peut conclure par le lemme \ref{recoll} comme \`a la section
suivante. Voici toutefois un argument alternatif.

Le cas particulier ${\cal F}= \mu_*\pi^* O_{\hat S}(-A)$, avec  $\mu: \hat S \to S$ d\'esingularisation,   $\hat \Si = \hat S \times _S \Si$, $A$  diviseur ample sur $\hat \Si$, donne lieu \`a:
\begin{coro} \label{lisbig}
Si $n\gg 0$,
$H^0_{2} (\hat \Si, \mu^*{\cal L}^n \otimes \pi^*O_{\hat S}(-A))\not= 0$.

En particulier,  $\mu^*{\mathcal L}$ poss\'ede une m\'etrique p\'eriodique de courbure
sup\'erieure \`a une classe K\"ahlerienne au sens des courants.
\end{coro}

\begin{coro}
On suppose de plus que pour tout sous espace propre $Z$ de S, ${\mathcal L}_{Z\times _S \Si}$ est $G$-nef. Alors $\mu^*\mathcal L$
est $G$-nef.  En particulier, pour tout $N>0$ il existe un faisceau d'ideaux de Nadel $I_N$ et $n_N$
sur
$\hat S$ tel que,  quelque soit $\mathcal M$ $G$-nef,
$n\ge n_N$ et $q\ge 1$
$$H^q_2 ( \hat \Sigma, K_{\hat \Sigma} \otimes
\mu^*({\mathcal L}^n (-NH)\otimes \mathcal M)\otimes \pi^*I_N)=0$$
\end{coro}
\begin{prv}Le premier point est cons\'equence de  \cite{P}, Th\'eor\`eme 1.C.3.
Le deuxi\`eme point r\'esulte de la version $L_2$ du th\'eor\`eme d'annulation de
Kawamata-Viehweg (voir  l'exploitation que  \cite{E1} fait de \cite{D1}).
\end{prv}

\begin{prop}\label{lisample}
Soit ${\mathcal L} \in Ob (L_G(\Sigma))$, tel que, pour tout
sous espace analytique de $S$, ${\mathcal L}^{\dim Z}.Z>0$.
Alors $\mathcal L$ est $G$-ample.
\end{prop}
\begin{prv}
Par r\'ecurrence sur $d=\dim \Si$, on peut supposer que pour tout sous espace propre $Z$ de S, ${\mathcal L}_{Z\times _S \Si}$ est $G$-ample. Alors pour tout $G,\chi$-fac de torsion
$T$, il existe $n_T$ tel que si $n\ge n_T$ et $q\ge 1$ et $\mathcal M$ est $G$-nef,
$H^q_2(\Si, T\otimes \mathcal M\otimes{\mathcal L}^n)=0$.

Le conoyau de $\mu_*(K_{\hat \Si}\otimes \pi^*I_N) \to \mu_*K_{\hat \Sigma}$ est un faisceau de torsion.
 Donc, pour tout $N\ge0$, il existe $n_N$ tel que, si $n\ge n_N$, $ q\ge 1$,
et $\mathcal M$ $G$-nef:
$$H^q_2(\Sigma, \pi^* (\mu_*K_{\Sigma}\otimes O(-NH))\otimes {\mathcal L}^n)=0
$$

Soit $\cal F$ un objet de $C_{G,\chi}(\Sigma)$. On peut trouver une r\'esolution de la forme:
$$
0 \to R^{-d}\to R^{-d+1} \to \ldots \to R^0 \to {\mathcal F} \to 0
$$
Pour $0\le i\le d-1$, $R^i$ est somme directe d'un nombre fini de faisceaux de la forme
$\pi^*(\mu_*K_{\hat \Sigma}\otimes\pi^* O_{S} (-NH) )\otimes {\mathcal L}^m$.

Pour calculer la cohomologie $L_2$ de $\mathcal F\otimes {\mathcal M}\otimes {\mathcal L}^n$,
on peut utiliser la suite spectrale d'hypercohomologie associ\'e \`a la r\'esolution
pr\'ec\'edente. Cette suite spectrale aboutit \`a
$H^*(\Si, \mathcal F \otimes \mathcal M \otimes \mathcal L^n)$ et son terme $E_2$
v\'erifie:
 $$E_2^{p,q}= H^p_2(\Si, R^q \otimes \mathcal M \otimes \mathcal L^n)$$

Or le point pr\'ec\'edent assure que, si $n\gg 0$,  $E_2^{p,q}=0$ si $p\not=-d$ et $q\not=0$
ou si $p=-d $ et $q\not\in[0,d]$. En particulier $E_2^{p,q}=0$ si $p+q\ge 1$ et, pour $n\gg0$
et $k\ge 1$, il suit que  $H^k(\Si, \mathcal F \otimes \mathcal M \otimes \mathcal L^n)=0$.
$\mathcal L$ est donc $G$-ample.
\end{prv}

\subsection{Cas o\`u $[\omega]\in H^2(\pi_1(X),\mathbb R)^{1,1}+ NS(X)\otimes {\mathbb R}$.}
\paragraph{}
On suppose que $\pi:\Sigma\to S$ est le rev\^etement universel d'un espace complexe r\'eduit projectif alg\'ebrique
compact $S$ de dimension $d$.
Soit $[\omega]\in H^1(S,{\mathcal H}_{\mathbb R})$ telle que:
\begin{itemize}
\item $[\omega]^{\dim Z}.Z>0$ pour $Z$ un sous espace alg\'ebrique de $S$.
\item Il existe
une suite d'objets de $L_G(\Si)$ $({\mathcal L}_k )_{k\in\mathbb N}$  et une suite d'entiers positifs
$(n_k)_k$ tels que $\lim_{k\to\infty}[c_1({\mathcal L}_k,h_k)]/n_k =[\omega]$.
\end{itemize}

\begin{lem}\label{petitedif} Soit $\omega_S$ une forme de K\"ahler sur
$S$.

Il existe une suite de r\'eels positifs tendant vers z\'ero $(\delta_k)_{k\in \mathbb N}$
et un repr\'esentant lisse $\Delta_k$ de $[\omega]-[c_1({\mathcal L}_k)]$
avec $-\delta_k \omega_S \le \Delta_k \le \delta_k \omega_S$.
\end{lem}

\begin{prop} \label{propo17}
$[\omega]$ est une classe de K\"ahler.
\end{prop}
\begin{prv} Par r\'ecurrence sur $d=\dim S$,
$[\omega]|_Z$ peut \^etre suppos\'ee de K\"ahler pour tout sous espace propre de $S$.

Fixons $H\subset S$ une section hyperplane. On peut trouver $k_0$ tel que, si
$k\ge k_0$ $\mathcal L_k |_H$ est $G$-ample et $\mathcal L_k^{\dim S}>0$ .

On fixe $O_{ S}(A)$ un diviseur ample
sur $S$ muni d'une m\'etrique lisse dont la courbure est la forme de K\"ahler $\omega_{ S}$.

Par le corollaire \ref{lisbig}, il existe $\epsilon \in \mathbb Q^+$ tel que pour $k\ge k_0$
on peut trouver une m\'etrique singuli\`ere \`a singularit\'es logarithmiques
$h_k$  sur ${\mathcal L}_k$ telle que
$\Theta(\mu^*L_k,h_k)/n_k\ge\epsilon \omega_{ S}$ au sens des courants.

Combinant ceci avec le lemme \ref{petitedif}, on trouve un repr\'esentant
big \`a singularit\'es alg\'ebriques de $[\omega]$. Invoquer le lemme \ref{recoll}
termine la preuve.
\end{prv}

\subsection{Fin de la preuve du Th\'eor\`eme \ref{theo1}}
\paragraph{}
On reprend les notations de l'introduction.
Soit $\xi:\Xi \to \Delta$ une d\'eformation de $X$.

Soit $O_{\Xi}(A)$ un faisceau inversible $\xi$-ample et $h$ une m\'etrique
 hermitienne lisse sur $O_{\Xi}(A)$ telle que $C_1(O_{\Xi}(A), h)$
 est lisse et strictement  positive sur chaque fibre $\Xi_t$ pour $t\in \Delta$
 assez petit.

Soit $(\omega_t)_{t\in \Delta}$ une famille continue de
$(1,1)$-formes ferm\'ees sur $\Xi_t$  avec $[\omega_{\Xi_0}]=[\omega]$
et $(t_n)_n$ une suite convergeant vers $0$ telle que:

$$[\omega_{t_n}]\in H^2(\pi_1(\Xi_{t_n}) ,\mathbb R)^{1,1} +NS(\Xi_{t_n})\otimes {\mathbb Q}$$

On suppose que pour tout sous espace analytique ferm\'e $Z\subset X$
$\omega^d. Z>0$.

Soit $\Upsilon \to \Delta$ une section hyperplane lisse de $\Xi \to
\Delta$.

Par r\'ecurrence sur $\dim_{\mathbb C} X$, $[\omega]_{\Upsilon_0}$ est une classe de K\"ahler sur $\Upsilon_0$.
En particulier, il existe $\epsilon \in \mathbb Q_{ >0}$ tel que
$([\omega] - \epsilon c_1(O_X(A)))|_{\Upsilon _0}$ est une classe de K\"ahler.

La stricte positivit\'e \'etant une condition ouverte, il suit qu'existe une fonction
$\phi$ sur $\Xi$, de classe
$C^{\infty}$, telle que
pour $t\in \Delta$ assez petit
$\omega_t -\epsilon C_1(O_{\Xi_t}(A), h) +dd^c\phi_{\Xi_t}$ est une
classe de K\"ahler sur $\Upsilon_t$.

Quitte \`a prendre $\epsilon$ plus petit,  le corollaire \ref{lisbig}
implique que $[\omega_{t_n}]-\epsilon c_1 (O_{\Xi_{t_n}}(A))$ est
repr\'esent\'ee par un courant positif ferm\'e d\'efini sur $\Xi_{t_n}$.

Par une variante ais\'ee du th\'eor\`eme de compacit\'e de Bishop,
 on peut passer \`a la limite pour obtenir que $\omega - \epsilon c_1(O_{X}(A))$
est repr\'esent\'ee par un courant positif ferm\'e.

Le th\'eor\`eme de r\'egularisation de Demailly \cite{D3} implique que
$\omega$ est repr\'esent\'ee par un courant positif ferm\'e $\ge \frac{\epsilon}2
c_1(O_X(A), h_X)$ lisse hors d'un ensemble analytique propre. Invoquer le
lemme \ref{recoll} termine la preuve.

Philippe Eyssidieux.

CNRS-UMR 5580. Laboratoire Emile Picard.

Universit\'e Paul Sabatier, UFR MIG.

118, Route de Narbonne.

31062 Toulouse Cedex (France)

\begin{verbatim}
 e-mail: eyssi@picard.ups-tlse.fr
\end{verbatim}

\end{document}